\documentclass[a4paper,notitlepage,11pt]{amsart}%
\usepackage{amssymb}
\usepackage{graphicx, enumerate}
\usepackage[shortlabels]{enumitem}
\usepackage{amsmath}
\usepackage{amsthm}
\usepackage{amsfonts}%
\usepackage{tikz}
\usepackage{dsfont}
\usepackage{mathrsfs}
\usepackage{color, soul}
\usepackage{multicol}
\usepackage{tcolorbox}
\usepackage{fancybox}

\usepackage{geometry}
\usepackage{etoolbox}
\let\svlabel\label
\let\svref\ref
\def\unusedlabels{}
\renewcommand\label[1]{\svlabel{#1}\global\edef\unusedlabels{\unusedlabels$<$#1$>$ }}
\renewcommand\ref[1]{\svref{#1}%
  \edef\teststring{$<$#1$>$}%
  \edef\tmp{\unusedlabels}%
  \def\unusedlabels{}%
  \expandafter\refhelper\tmp\relax%
}
\def\refhelper#1 #2\relax{%
  \edef\expandcase{#1}%
  \ifdefstrequal{\teststring}{\expandcase}{}{\edef\unusedlabels{\unusedlabels#1 }}%
  \if\relax#2\relax\else\refhelper#2\relax\fi
}

\usepackage[colorlinks=true,linkcolor=blue,urlcolor=blue,citecolor=blue]{
hyperref}

\def\CC{\mathbb{C}}

\def\KK{\mathbb{K}}

\def\NN{\mathbb{N}}

\def\XX{\mathbb{X}}

\def\mcB{\mathcal{B}}
\def\mcE{\mathcal{E}}

\def\mcH{\mathcal{H}}

\def\1{\mathds{1}}

\def\bC{\mathbf{C}}

\def\xx{\mathbf{x}}

\providecommand{\U}[1]{\protect\rule{.1in}{.1in}}

\theoremstyle{plain}
\newtheorem{theorem}{Theorem}[section]

\newtheorem{corollary}[theorem]{Corollary}

\newtheorem{proposition}[theorem]{Proposition}

\theoremstyle{definition}
\newtheorem{definition}[theorem]{Definition}
\newtheorem{remark}[theorem]{Remark}
\numberwithin{equation}{section}
\numberwithin{theorem}{section}

\usepackage{enumitem}

\newcommand{\abs}[1]{\left| #1 \right|}
\newcommand{\absito}[1]{\big| #1 \big|}

\newcommand{\normvar}[1]{\Big\Vert #1 \Big\Vert}
\newcommand{\normita}[1]{\big\Vert #1 \big\Vert}

\def\vareptil{\tilde\varepsilon}

\def\supp{{\rm supp}}

 \thanks{This work was supported in part by ANPCyT PICT 2018-04104. The first author's research was also partially supported by CONICET PIP 11220200102366CO and UBACyT 20020220300242BA and the second author by CONICET PIP 11220200101609CO and PAI-UdeSA 2023.}


\begin{document}
 	
	\title[The fundamental functions of the canonical basis  $\{e^{-\lambda_ns}\}$ of $ \mcH_p^\lambda$]{The fundamental functions of the canonical basis \\ of Hardy spaces of Dirichlet series}
	\date{}

\author[D. Carando]{Daniel Carando}
\address{Daniel Carando\\
Departamento de Matem\'atica 
\\Facultad de Ciencias Exactas y Naturales\\ Universidad de Buenos Aires 
and IMAS (UBA--CONICET) 
\\ Pabellón I, Ciudad Universitaria \\
(1428) Buenos Aires,   Argentina}
\email{dcarando@dm.uba.ar}

\author[S. Lassalle]{Silvia Lassalle}
\address{Silvia Lassalle\\
Departamento de Matem\'atica\\
Universidad de San Andr\'es and
IMAS--CONICET\\ Vito Duma 284 (1644) Victoria, Buenos Aires, Argentina}
\email{slassalle@udesa.edu.ar}

\author[L. Milne]{Leandro Milne}
\address{Leandro Milne\\
Departamento de Matem\'atica 
\\Facultad de Ciencias Exactas y Naturales\\ Universidad de Buenos Aires 
\\ Pabellón I, Ciudad Universitaria \\
(1428) Buenos Aires,  Argentina}
\email{lmilne@dm.uba.ar}
 
 \maketitle

\begin{abstract}
Given a frequency $\lambda=(\lambda_n)$, we consider the Hardy spaces $ \mcH_p^\lambda$ of $\lambda$-Dirichlet series \( D = \sum_n a_n e^{-\lambda_n s}\) and study the asymptotic behavior of the upper and lower democracy functions of its canonical basis $\mathcal B=\{e^{-\lambda_ns}\}$.
For the ordinary case, $\mathcal B=\{n^{-s}\}$, we give the correct asymptotic behavior of all such functions, while in the general case we give sharp lower and upper bounds for all possible behaviors.  Moreover, for $p>2$  we present examples showing that any intermediate behavior (between the extreme bounds) can occur. We also study how different properties of the frequency $\lambda$ lead to particular behaviors of the corresponding fundamental functions. Finally, we apply our results to analyze greedy-type properties of $\mathcal B=\{e^{-\lambda_ns}\}$ for some particular $\lambda$'s.
\end{abstract}

\section*{Introduction}

The study of general Dirichlet series from a more analytic point of view  dates back to the beginning of the 20th century with its groundbreaking contribution to analytic number theory. The general theory was developed by Hadamard, Landau, Hardy, Riesz and Bohr, among others, where the notion of series given by a frequency $\lambda$ even appears. It can be said that the topic experienced a revival with the seminal paper by Hedenmalm, Lindqvist, and Seip~\cite{hedenmalm1997}, where the authors provided new insights by introducing the first \emph{Hardy} spaces of Dirichlet series to solve a problen in harmonic analysis. 
In the last two decades, the topic has attracted the attention of several authors leading to the development of the theory of Hardy spaces of Dirichlet series, first with the work of Bayart (ordinary series)~\cite{Bayart02} and later with the work of Defant and Schoolmann (general series)~\cite{DeSch}. 

This paper investigates the fundamental functions associated with the canonical basis $\mathcal B=\{e^{-\lambda_ns}\}$ of the different Hardy spaces of general Dirichlet series, determined by a frequency $\lambda=(\lambda_n)$. These functions, namely the lower/upper democracy and super-democracy functions, play an important role in understanding different properties of a basis.
We begin by focusing on the classical Hardy spaces $\mathcal{H}_p$  of ordinary Dirichlet series, which is done in Section~\ref{sec-ordinarias}. Theorem~\ref{teo-democ-ordinarias} establishes the precise asymptotic behavior of all fundamental functions in this setting. The estimations obtained, reveal that the canonical basis exhibits ``good properties" within $\mathcal{H}_p$  only for $p=2$, i.e. when the basis is actually orthonormal, see Corollary~\ref{cor-supercoro}.

Section~\ref{sec-generales} extends the analysis to Hardy spaces $\mcH_p^\lambda$ of  $\lambda$-Dirichlet series for some frequency $\lambda=(\lambda_n)$. These are series of the form 
$$ 
D = \sum_n a_n e^{-\lambda_n s}.
$$
As expected, the asymptotic behaviors of the fundamental functions become dependent on the chosen frequency. Theorem~\ref{teo-fi-series-generales} provides sharp lower and upper bounds for all possible behaviors. We also show explicit frequency examples that achieve the extreme behaviors. 
Furthermore, for $p>2$, we adapt a construction from Bourgain's work~\cite{Bo} (see Proposition~\ref{prop-Bourgain}) to demonstrate that any intermediate behavior can be realized. Also, for specific frequency families exhibiting properties like containing long arithmetic progressions, being linearly independent over $\mathbb Q$, or being lacunary, the precise asymptotic behavior can be determined.

Section~\ref{sec-unconditional} delves into unconditional and democratic properties of the canonical basis. In the general setting, ``good properties" within a particular $\mcH_p^\lambda$ space do not necessarily imply $p=2$. Instead, they may indicate an isomorphism between $\mcH_p^\lambda$ and $\mcH_2^\lambda$ (see Proposition~\ref{prop-uncond}) or a close relationship between them (see Propositions~\ref{prop-democ-aprox-h2} and Theorem~\ref{teo-qglc_dn}).

Theorem~\ref{teo-qglc_dn} states, in particular, that if the canonical basis in $\mcH_p^\lambda$ is suppression unconditional for constant coefficients (SUCC), then it is democratic, an implication that is not generally true (see Section~\ref{sec-prelim} for the definitions). Considering the key role that democracy plays in the interplay of the different types of greedy bases, in Section~\ref{sec-greedy}, we employ the estimates obtained for (super-)democratic functions and apply them to this field. For frequencies containing arbitrarily long arithmetic progressions, we prove that $\mcB$  is SUCC in $\mcH^\lambda_p$ if and only if $p=2$. This leads us to conclude that for such frequencies (in particular, for ordinary Dirichlet series and Fourier series) $\mcB$ does not enjoy any of the most commonly used greedy properties unless we stay in the Hilbert space $\mcH^\lambda_2$. 

\medskip
In Section~\ref{sec-prelim} we give all the necessary definitions and state many of the general results used in the sequel. We refer the reader to~\cite{defant2019dirichlet} and~\cite{DeSch} for more information on Hardy spaces of  Dirichlet series and to~\cite{temlyakov1998greedy} and~\cite{albiac2021dissertationes} for further reading on fundamental functions and greedy-type bases.

\section{Definitions and preliminary results} 
\label{sec-prelim}

The notions of democracy and the associated fundamental functions are defined 
with respect to a coordinate system that we will call \emph{basis}. These concepts, which are closely related to different notions of greedy algorithms, only appeal to a fixed finite number of elements of the system. 

Let $\XX$ be an infinite dimensional separable Banach space over the real or complex field $\KK$.
A sequence $\mcB=\{\xx_n\}$ is a \emph{Markushevich basis} for $\XX$ if it spans a dense subspace of $\XX$
and admits a (unique) \emph{total} sequence of \emph{biorthogonal functions} $\mcB^*=\{\xx_n^*\}$ in the dual space $\XX^*$. That is,  
$\xx_n^*(\xx_j)=\delta_{n,j}$ for all $n,j\in \NN$ and, whenever
$\xx_n^*\left(f\right)=0$ for\ all $n\in \NN$, we necessarily have $f=0$.
As usual, we refer to $\mcB^*$ as the \emph{dual basis} of $\mcB$. Also, we use $\supp(f)$ to denote the support of $f\in \XX$, that is the set $\{j\in \NN: \xx^*_j(f)\ne 0\}$.

When there is $\bC>0$ such that for all $f\in \XX$ and all $n\in \NN$
$$
\normita{\sum_{j=1}^{n}\xx_i^{*}(f)\xx_i}\le \bC\|f\|,
$$
the sequence $\mcB$ is a \emph{Schauder basis} and the minimum $\bC>0$ for which the above inequality holds is the \emph{basis constant} of $\mcB$. A \emph{basic sequence} is a Schauder basis for the closure of its span. A Schauder basis $\mcB$ is $\bC$-\emph{unconditional} if for all $x\in \XX$, all  $(a_j)\subset\KK\colon\ |a_j|\le 1$ for all $j\in \NN$, 
$$
\normita{\sum_{j=1}^{\infty}a_j\xx_j^*(x)\xx_j}\le \bC\|x\|.
$$

Somehow related to the notion of unconditionality we have the concepts of (super-)democracy. For the definitions, we need the following notation. 

Given any set $A\subset \NN$, we define
$$
\mcE_A:=\{\varepsilon=(\varepsilon_j)_{j\in A}\colon\ |\varepsilon_j|=1\quad \forall j\in A\}
$$
with the convention that $\mcE_A=\emptyset$ if $A=\emptyset$. For a basis $\mcB=\{\xx_n\}$, given a finite set 
$A\subset \NN$ and $\varepsilon\in\mcE_A$, we denote 
$$
\1_{\varepsilon, A}:= \1_{\varepsilon, A,}[\mcB,\XX]=\sum_{j\in A}\varepsilon_j\xx_j, 
$$
with the convention that any sum over the empty set is zero. Also, if $\varepsilon\in \mcE_A$ and $B\subset A$, we write $\1_{\varepsilon,B}$ considering the natural restriction of $\varepsilon$ to $B$.  If $\varepsilon_j=1$ for all $j$, we simply write $\1_A$. 

\begin{definition} Let $\mcB$ be a basis for $\XX$ and $\bC>0$.  Then: 
\begin{enumerate}[(a)]
\item $\mcB$ is \emph{democratic with constant $\bC$} (or $\bC$-democratic) if for any finite sets $A, B \subset \NN$ with $|A| \le |B|$ the following holds
\begin{align}\label{def-demo}
\|\1_{A}\|\le \bC \|\1_{B}\|.
\end{align}
\item $\mcB$ is \emph{super-democratic with constant $\bC$} (or $\bC$-super-democratic)  if for any finite sets $A, B \subset \NN$ with $|A| \le |B|$, any $\varepsilon\in \mcE_A $ and $\vareptil\in \mcE_B$ the following holds
\begin{align}\label{def-sup-demo}
\|\1_{\varepsilon,A}\|\le \bC \|\1_{\vareptil,B}\|.
\end{align}
\item $\mcB$ is \emph{supression unconditional for constant coefficients with constant} $\bC$ (or $\bC$-SUCC) if \eqref{def-sup-demo} holds for any finite sets $A \subset B \subset \NN$. 
\end{enumerate}
\end{definition}

Let us recall the interplay among these properties. Unconditionality and democracy are two concepts that do not imply one another (see, e.g.,  \cite[Sec~1.3.~p.23]{temlyakov2011greedy}). The same holds for SUCC and democracy properties. Indeed, since any unconditional basis is SUCC, the example given in \cite{temlyakov2011greedy} shows one assertion. A democratic basis that fails to be SUCC was given, for instance, in \cite[Example~11.19]{albiac2021dissertationes}. On the other hand, it follows from the definitions that any super-democratic basis is democratic and SUCC. In addition, a SUCC and democratic basis is super-democratic, see for instance \cite[Proposition~5.1]{albiac2021dissertationes}. 

Now, we introduce the upper and lower democracy functions \(\varphi_{u}\) and \(\varphi_{l}\) respectively, which are useful to study how close to democracy a basis $\mcB$ is. For \( N \in \mathbb{N} \), we define its \emph{fundamental functions} as
$$
\begin{array}{rcl}
    \varphi_{u}(N)=&\varphi_{u}[\mathcal B, \mathbb X](N):=&\sup\{\normita{\mathds{1}_A}\colon A\subset \NN , \abs{A}\le N\}, \\
    \varphi_{l}(N)= &\varphi_{l}[\mathcal B, \mathbb X](N):=&\inf\{\normita{\mathds{1}_A}\colon A\subset \NN , \abs{A}\ge N\}. 
\end{array}
$$
We also consider the  upper and lower super-democracy functions of \(\mcB\) defined, respectively, for \( N \in \mathbb{N} \) by
$$\begin{array}{rcl}
    \varphi_{u,\mcE}(N) = &\varphi_{u,\mcE}[\mathcal B, \mathbb X](N):= & \sup\{\normita{\mathds{1}_{\varepsilon,A}}\colon A\subset \NN, \abs{A}\le N, \varepsilon\in \mcE_A\}, \\
    \varphi_{l,\mcE}(N)= &\varphi_{l,\mcE}[\mathcal B, \mathbb X](N):= & \inf\{\normita{\mathds{1}_{\varepsilon,A}}\colon A\subset \NN , \abs{A}\ge N, \varepsilon \in \mcE_A\}. 
\end{array}
$$
The functions  \(\varphi_{u},\, \varphi_{l}, \, \varphi_{u,\mcE},\, \varphi_{l,\mcE} \) are non-decreasing and satisfy
\begin{equation*}
    \varphi_{l,\mcE}\le \varphi_{l}\le \varphi_{u}\le \varphi_{u,\mcE}.
\end{equation*}
For Banach spaces the functions \( \varphi_{u} \)  and \(\varphi_{u,\mcE} \) are equivalent, see equation (8.3) in~\cite{albiac2021dissertationes}, and for the class of SUCC bases the equivalence reaches the functions \( \varphi_{l} \) and \(\varphi_{l,\mcE} \). 

\subsection{Hardy spaces of Dirichlet series}

A \emph{Dirichlet series} is a (formal) sum of the type \[ f(s) = \sum_{n=1}^{\infty} a_n n^{-s}, \]
where \( (a_n)\) is a sequence of complex coefficients and \( s \) is a complex variable. Dirichlet series and power series in infinitely many variables can be identified via the Bohr transform as follows.

Denote by \( \mathbb{Z}^{(\mathbb{N})} \)  the set of all sequences \( \alpha = (\alpha_1, \alpha_2, \dots )\) such that all entries are integers and only a finite number of them are nonzero, and denote by  $\mathbb{N}_0^{(\mathbb{N})}$ the subset of \( \mathbb{Z}^{(\mathbb{N})} \) consisting of the elements of non-negative entries.
For \( z\in \mathbb{C}^\infty \) and a multi-index \( \alpha \in \mathbb{Z}^{(\mathbb{N})} \), we write \( z^\alpha = z_1^{\alpha_1}z_2^{\alpha_2}\dots z_n^{\alpha_n} \), where \( n= \max \{j:\alpha_j\neq0\} \).
A formal power series in infinitely many variables is an expression of the form  \[ \sum_{\alpha\in \mathbb{N}_0^{(\mathbb{N})}}   c_\alpha z^\alpha. \]

We next identify the space of all formal power series \( \mathfrak{P} \) with the space of all formal Dirichlet series \( \mathfrak{D} \). Let \( \mathfrak{p}=(p_k)_k \) be the sequence of prime numbers. For each \( n\in \NN\) we take its prime number decomposition  \( n = p_1^{\alpha_1} \dots p_k^{\alpha_k}= \mathfrak{p}^\alpha\), with $\alpha=(\alpha_1,\dots,  \alpha_k,0,\dots)\in \mathbb{N}_0^{(\mathbb{N})}.$  Now, the identification is given by the bijective mapping
 \begin{center}
	\begin{tikzpicture}[node distance=0.1cm]
		\node (P) {\( \mathfrak{P} \)};
		\node (D) [,right of=P, xshift=4cm] {\( \mathfrak{D} \)};
		\node (empty) [above of=P, yshift=0.1cm] {};
		\node (empty2) [below of=P, yshift=-0.4cm] {};
		\node (B) [right of=empty, xshift=2cm] {$\mathfrak{B}$};
		\node (coef) [right of=empty2, xshift=2cm] {\footnotesize{\(c_\alpha=a_{\mathfrak{p}^\alpha}\)}};
		\node (powerserie) [below of=P, yshift=-0.5cm] { \(\sum_{\alpha} c_\alpha z^\alpha\)};
		\node (dirichletserie) [right of=powerserie, xshift=4cm] {\(\sum_n a_n n^{-s}\)};
		\draw [->] (P) -- (D);
		\draw [->] (powerserie) -- (dirichletserie);
	\end{tikzpicture}
\end{center}
which is known as the \emph{Bohr transform}.
	
We consider the infinite dimensional poly-torus $\mathbb{T}^{\infty} = \prod_{k=1}^\infty \mathbb{T}$ with its normalized Lebesgue measure, which is the countable product measure  of the normalized Lebesgue measure
on $\mathbb{T}$. For any  multi-index $\alpha = (\alpha_1, \dots, \alpha_n,0, \ldots ) \in \mathbb{Z}^{(\mathbb{N})}$ 
the $\alpha$th Fourier coefficient $\hat{f}(\alpha)$ of $f\in L_{1}(\mathbb{T}^{\infty})$ is given by
\[
\hat{f}(\alpha) = \int_{\mathbb{T}^{\infty}} f(z) z^{- \alpha} dz\,.
\]

The \emph{Hardy space} $H_p(\mathbb{T}^\infty)$, $1\le p\le \infty$, is the set of all functions  \( f \in L_p(\mathbb{T}^\infty) \) such that  \( \hat{f}(\alpha) = 0 \) for all \( \alpha \in \mathbb{Z}^{(\mathbb{N})} \setminus \mathbb{N}_0^{(\mathbb{N})}  \).
Each \( f \in H_p(\mathbb{T}^\infty) \) is uniquely determined by its Fourier coefficients $\{\hat{f}(\alpha)\}_{\alpha}$  and thus it is associated to a formal power series:  
\[ 
f(z) \sim \sum_{\alpha \in \mathbb{N}_0^{(\mathbb{N})} } \hat{f}(\alpha) z^{\alpha}. 
\] 
In other words, \( H_p(\mathbb{T}^\infty) \) can be viewed as a subspace of \( \mathfrak{P}. \) 
  
For \( 1\le p \le \infty \), the \emph{Hardy space of Dirichlet series}  \( \mathcal{H}_p, \) is defined  as the image of  \( H_p(\mathbb{T}^\infty) \) by the Bohr transform \( \mathfrak{B} \):
\begin{equation}
		\mathcal{H}_p := \mathfrak{B}(H_p(\mathbb{T}^\infty)).
	\end{equation} 
On \( \mathcal{H}_p \) we consider the norm defined by    
\begin{equation*}
        \normita{D}_{\mathcal{H}_p} := \normita{\mathfrak{B}^{-1}(D)}_{H_p(\mathbb{T}^\infty)}, \ \text{ for all}\ D\in \mathcal{H}_p. 
    \end{equation*}
Thus, the Bohr transform $\mathfrak{B}$ is an isometric isomorphism between $H_p(\mathbb{T}^\infty)$ and $\mathcal{H}_p$. A \emph{Dirichlet
polynomial} $D\in \mathcal{H}_p$ is a Dirichlet series with only finite non-zero coefficients, say $D = \sum_{n\in A}
a_n n^{-it}$, for $A\subset \mathbb N$ a finite set and $(a_n)_{n\in A}$ complex numbers. In this case, we have an intrinsic way to
compute the norm of $D$:
\begin{equation}\label{eq-norma-dir-pol}
        \normita{D}_{\mathcal{H}_p}= \left(\lim_{R\to \infty} \frac{1}{2R} \int_{-R}^{R} 
        \Big|\sum_{n\in A} a_n n^{-it}\Big|^pdt\right)^{\frac{1}{p}} \quad \text{and} \quad \normita{D}_{\mathcal{H}_\infty}= \sup_{t\in \mathbb R } \Big|\sum_{n\in A} a_n n^{-it}\Big| 
    \end{equation} 
This gives, for $1\le p<\infty$, an alternative definition of $\mathcal{H}_p$ as the completion of the space of Dirichlet polynomials with the norm given in \eqref{eq-norma-dir-pol}, see~\cite{Bayart02} or~\cite{defant2019dirichlet}. In the case $p=\infty$, the completion of the space of Dirichlet polynomials is a closed subspace of $\mcH_\infty$ that we denote by  $\mathcal{A}_\infty$. For a characterization of this space see Proposition 2.2 and Theorem 2.3 of~\cite{aron-bayart-Gauthier-Maestre-Nestoridis-2017dirichlet}, where the notation  $\mathcal{A}(\CC_+)$ is used.

The sequence \( \mcB := \{  n^{-s}\}\) is a Schauder basis  of \( \mathcal{H}_p \), \( 1 < p < \infty\), see~\cite{aleman-Olsen-Saksman-2014fourier}, and it is a Markushevich basis for $\mcH_1$ and  \( \mathcal{A}_\infty\), for \( p=1\) and \( p = \infty\), respectively. 
Indeed, by definition we have $\mathcal{A}_\infty=\overline{[n^{-s}: n\in\NN]}$ and $\mcH_1=\overline{[n^{-s}: n\in\NN]}$ (see for instance \cite[Theorem~1.10]{defant2019dirichlet}). In both cases, the biorthogonal functions are the coordinate functions. The rest of the conditions are straightforwardly satisfied. In what follows, we refer to \( \mcB = \{  n^{-s}\}\) as the \emph{canonical basis} of the corresponding space.  
	 
\begin{remark} The Fourier basis is unconditional in $L_p(\mathbb T)$ if and only if $p=2$ (see, for example, \cite[IID8~Prop.~9]{LibroWoj}. As a consequence, the basis \( \mathcal B \) is unconditional in $\mathcal{H}_p$ if and only if $p= 2 $ (see also  \cite[Remark~4]{temlyakov1998greedy} or \cite[Proposition~4]{CaDeSe18}).
\end{remark}

\subsection{Some useful inequalities}

In this section, we recall some inequalities that will be useful in the sequel. As usual, for \( 1\le r \le \infty \) we denote by \( r'\) the conjugate index of $r$, i.e., \( \frac{1}{r}+ \frac{1}{r'}=1 \) with the usual convention that $r'=1$ if $r=\infty$. 

The classical Hausdorff--Young inequalities can be extended to the Dirichlet series framework (see, for example,~\cite{carando2020hausdorff, GarciaCuervaetal}). Thus, for every \( 1\le p \le 2, \) and all scalars $(a_n)$, we have
\begin{equation}\label{hausdorff-young:1-2}
\normvar{\sum_{n=1}^{\infty} a_nn^{-s}}_{\mathcal{H}_{p'}} \le \Big( \sum_{n=1}^{\infty} \abs{a_n}^p \Big)^{\frac{1}{p}},
\end{equation}
and for \( 2\le q \le \infty \),  
\begin{equation}\label{hausdorff-young:2-infty}
\Big(\sum_{n=1}^{\infty} \abs{a_n}^q \Big)^{\frac{1}{q}} \le \normvar{\sum_{n=1}^{\infty} a_nn^{-s}}_{\mathcal{H}_{q'}}. 
\end{equation}

Let us state Khinchin inequalities in their probabilistic version. For this, let $\{\delta_n\}$  be a sequence of independent Bernoulli random variables, each taking the values $-1$ or $1$ with probability $1/2$. Then, given $1\le  p < \infty$ there exist (best) constants $\mathbf{k}_p$ and $\mathbf{K}_p$ such that
\begin{equation}\label{eq-kinchin}
\mathbf{k}_p^{-1} \left(\sum_{n\in A}  \left\vert a_n\right\vert^2 \right)^\frac{1}{2}\le \left[\mathbb{E}\left(\Big| \sum_{n\in A}  a_n \delta_n \Big|^p \right)\right]^\frac{1}{p}\le \mathbf{K}_p\left( \sum_{n\in A} \vert a_n \vert^2 \right)^\frac{1}{2},
\end{equation} 
for any finite set $A\subset \NN$ and all $(a_n)_{n\in A} \subset \mathbb C.$

We also recall the Khinchin--Steinhaus inequalities (see Theorem 6.8 in~\cite{defant2019dirichlet}). 
For each \( 1\le p < \infty, \) there are (best) constants \( \mathbf{s}_p, \mathbf{S}_p\ge 1 \) such that 
for any finite set $A\subset \NN$ and all $(a_n)_{n\in A} \subset \mathbb C$, we have 

\begin{equation}\label{eq-kin}
\mathbf{s}_p^{-1}\left(\sum_{n\in A}  \abs{a_n}^2\right)^{\frac{1}{2}} \le \left(\int_{\mathbb{T}^{\infty}} \Big|{\sum_{n\in A}  a_n z_n}\Big|^p dz\right)^{\frac{1}{p}}\le \mathbf{S}_p\left(\sum_{n\in A}  \abs{a_n}^2\right)^{\frac{1}{2}}.  \end{equation}	
The integral above can be interpreted as an expected value, so this inequality is a reformulation of \eqref{eq-kinchin} for Steinhaus random variables. 

Note that we can easily reformulate \eqref{eq-kin} in terms of Dirichlet series. 
Using the Bohr transform  and then the Khinchin--Steinhaus inequality, for $1\le p < \infty$  we have, for any finite set $A\subset \NN$ and all $(a_n)_{n\in A} \subset \mathbb C$,
\begin{equation}\label{sum first N primes}
\normvar{\sum_{n\in A} a_n p_n^{-s}}_{\mathcal{H}_p}= \normvar{\sum_{n\in A} a_n  z_n }_{H_p(\mathbb{T}^{\infty})} \approx \left(  \sum_{n\in A}  \abs{a_n }^2 \right)^{\frac{1}{2}} . 
\end{equation}
On the other hand, for $p=\infty$ we have $$\normvar{\displaystyle\sum_{n\in A} a_n p_n^{-s}}_{\mathcal{H}_\infty} =\sum_{n\in A}  |a_n|.$$ 

We also need estimates for the $p$-norm of the so-called Dirichlet kernel \( \mathbf{D}_N \), defined as	
\[
\mathbf{D}_N(z) = \sum_{j=1}^{N} z^j, \ \text{for } \ z \in \mathbb{\mathbb{T}}.
\]
Using \cite[Lemma~1.7]{duoandikoetxea2001fourier} for $p=1$ and \cite[Lemma~2.1]{anderson2007exponential} for $ p>1$, we have:
	\begin{equation}\label{eq-norm-dir-ker}
	    \normita{\mathbf{D}_N}_{L_p(\mathbb{T})} \approx \begin{cases}
		\log N & \text{if}\ \ p=1,\\
		  N^{\frac{1}{p'}} & \text{if}\ \ p>1.
	\end{cases}
	\end{equation} 

Again, these estimates  have easy  reformulations for Dirichlet series, now using only powers of 2:
		\begin{equation}\label{eq-norma-pot-de-2}
		     \normvar{\sum_{j=1}^{N}  (2^{j})^{-s} }_{\mathcal{H}_p} 
                = \left(\int_{\mathbb{T}}  \Big|\sum_{j=1}^{N} z_1^j \Big|^p dz_1 \right)^{1/p} 
                =  \normita{\mathbf{D}_N}_{L_p(\mathbb{T})} \approx \begin{cases}
		\log N &\text{if}\ \ p=1,\\
		  N^{\frac{1}{p'}} &\text{if}\ \ p>1.
	\end{cases}
		\end{equation}

Finally, we recall the definition of the Rudin--Shapiro polynomials (see, for instance, \cite[Proposition~9.7]{defant2019dirichlet}). 

\begin{definition}\label{rud-shap}
	The Rudin--Shapiro polynomials \( P_k\colon\mathbb{C}\rightarrow \mathbb{C} \) are defined recursively as follows: 
	\begin{align*}
		&	P_0(z):= z, \\ 
		&		P_k(z):= P_{k-1}(z^2)+\frac{1}{z} P_{k-1}(-z^2),\ \text{for} \ k \in \mathbb{N}.
	\end{align*}
\end{definition}	
By induction it can be easily checked that, for every \( k \) and \( z \),
\[ P_k(z)= \sum_{j=1}^{2^{k}} a_{kj}z^j, \]
where the coefficients \(  a_{kj} \) are either or \( 1 \) or \( -1 \). The important feature of these polynomials is that their supremum norm on $\mathbb T$  satisfy \( \normita{P_k}_{\infty} \le \sqrt{2^{k+1}} \), which are small compared with the norm of a general polynomial with coefficients of modulus 1.

\section{The canonical basis of \( \mathcal{H}_p \)} 
\label{sec-ordinarias}

The main goal of this section is to give asymptotically correct estimates of the democracy and super-democracy functions for the canonical basis $\mathcal B=\{n^{-s}\}$ of $\mathcal{H}_p$. We fix the following notation, for $1\le p \le \infty$ and $N\in \NN$,
\begin{align*}
\varphi_{u}^p(N):=\varphi_u[\mathcal B,\mathcal{H}_p](N),\ \,&\qquad \varphi_{l}^p(N):=\varphi_l[\mathcal B,\mathcal{H}_p](N), \\ \varphi_{u,\mcE}^p(N):=\varphi_{u,\mcE}[\mathcal B,\mathcal{H}_p](N),\ \,&\qquad \varphi_{l,\mcE}^p(N):=\varphi_{l,\mcE}[\mathcal B,\mathcal{H}_p](N).
\end{align*}

The following is the main result of this section.

\begin{theorem}\label{teo-democ-ordinarias} {The fundamental functions of the canonical basis $\mcB=\{n^{-s}\}$ in $\mathcal{H}_p$ satisfy:} 
\begin{enumerate}[{\rm (a)}]
\item If $p=1$,\vspace{-10pt}
\begin{eqnarray*} 
{\rm (i)}\quad  
\varphi_{l}^1(N) \approx \varphi_{l,\mcE}^1(N) \approx  \log(N), 
&\quad&  
{\rm (ii)}\quad 
\varphi_{u}^1(N) \approx \varphi_{u,\mcE}^1(N) \approx N^{\frac 12}. 
\end{eqnarray*}
\item If $1< p <\infty$,\vspace{-10pt}
\begin{eqnarray*}
{\rm (i)}\quad  
\varphi_{l}^p(N) \approx \varphi_{l,\mcE}^p(N) \approx  N^{\min \left\{\frac 1 2, \frac 1 {p'} \right\}}, 
&\quad& 
{\rm (ii)}\quad 
\varphi_{u}^p(N) \approx \varphi_{u,\mcE}^p(N) \approx  N^{\max \left\{\frac 1 2, 
       \frac 1 {p'} \right\}}. \label{eq-fip}
\end{eqnarray*}
\item If $p = \infty$,\vspace{-10pt}
\begin{eqnarray*} 
{\rm (i)}\quad  
\varphi_{l}^\infty (N) = N \quad \text{and}\quad \varphi_{l,\mcE}^\infty (N) \approx N^{\frac 1 2},  &\quad& 
{\rm (ii)}\quad  
\varphi_{u}^\infty (N)=\varphi_{u,\mcE}^\infty (N)=N. 
\end{eqnarray*}
\end{enumerate}
\end{theorem}

Before giving the proof of Theorem~\ref{teo-democ-ordinarias}, we present some simple consequences.
\begin{remark} \label{rem-conseq}
In light of the fact that a basis $\mathcal B$ is democratic if and only if the associated upper and lower democracy functions are asymptotically equivalent, Theorem~\ref{teo-democ-ordinarias} shows that the canonical basis of $\mathcal{H}_p$ is democratic if and only if $p=2$ or $p=\infty$. It is curious that for $p=\infty$ we get a \emph{good} behavior of the basis  despite the fact that, in this case, $\mathcal B$ even fails to be a basic sequence. About super-democracy, the theorem shows that it is only satisfied for $p=2$. We will see in Section~\ref{sec-unconditional} that, unless $p=2$, the basis $\mathcal B$ is not even suppression unconditional for constant coefficients. 
\end{remark}

\begin{proof}[Proof of Theorem~\ref{teo-democ-ordinarias}]
We start with (a)(\rm i). 
 From \cite[Corollary~1]{mcgehee1981hardy} 
   we have  
    \begin{equation}\label{inequality from mcgehee}
        \normvar{\sum_{k=1}^N c_k e^{in_kt}}_1\ge C^{-1} \log (N),
    \end{equation}
   whenever  \( \abs{c_k}\ge 1 \) for all \(k\), and \( n_1,\dots,n_N \) are all different integers. Given \(A \subset \mathbb{N}\) with \(\abs{A}\ge N \), via the Bohr transform, we have  
   \begin{equation*}
    \normvar{\sum_{n\in A} \varepsilon_n n^{-s}}_{\mathcal{H}_1} 
    = \normvar{\sum_{\alpha \in F} \varepsilon_\alpha z^{\alpha}}_{H_1(\mathbb{T}^\infty)} = \normvar{\sum_{\alpha \in F} \varepsilon_\alpha z^{\alpha}}_{H_1(\mathbb{T}^M)} 
    = \int_{\mathbb{T}^M} \Big|\sum_{\alpha \in F} \varepsilon_\alpha  z^{\alpha} \Big| dz
   \end{equation*}
   where  \( F= \{ \alpha \in \mathbb{N}^{(\mathbb{N})}\colon \mathfrak{p}^\alpha \in A \} \) and \( M = \max\{j\colon \alpha_j\ne 0, \alpha \in F\} \).
   Proceeding as in the proof of \cite[Proposition~2.4]{carando2016some}, see  (11) and (12) in~\cite{carando2016some},  we get 
    \begin{equation*}
          \int_{\mathbb{T}^M} \Big|\sum_{\alpha \in F} \varepsilon_\alpha z_1^{\alpha_1}\cdots z_M^{\alpha_M} \Big| dz 
        = \int_{\mathbb{T}^{M-1}} \left(\int_{\mathbb{T}}\Big|\sum_{\alpha \in F} \varepsilon_\alpha z_2^{\alpha_2}\cdots z_M^{\alpha_M} z_1^{\alpha_1+(m+1)\alpha_2+\cdots + (m+1)^{M-1}\alpha_M } \Big|  dz_1\right)dz_2\cdots dz_M 
         \end{equation*}
         where \(m= \max \{\alpha_j: \alpha\in F\} \). Since $\abs{\varepsilon_\alpha z_2^{\alpha_2}\cdots z_M^{\alpha_M} }=1$ for every $\alpha$ and the powers of \(z_1\) are all different, we can apply \eqref{inequality from mcgehee} to obtain
\begin{equation*}
    \int_{\mathbb{T}^M} \Big|\sum_{\alpha \in F} \varepsilon_\alpha z_1^{\alpha_1}\cdots z_M^{\alpha_M} \Big| dz \ge C^{-1} \log(|A|) \int_{\mathbb{T}^{M-1}} dz_2\cdots dz_M = C^{-1}\log(|A|) \ge C^{-1}\log(N).
\end{equation*}
Thus, \(\varphi_{l,\mcE}^p(N) \gtrsim \log(N).\) Moreover, by \eqref{eq-norma-pot-de-2} we have $\varphi_{l}^p(N) \lesssim \log N$, yielding the estimations in (a)(\rm i). 

\medskip
We handle together the upper functions in (a)(ii)  
and (b)(ii). 
By \eqref{sum first N primes} and \eqref{eq-norma-pot-de-2} we have $\varphi_{u}^p(N) \gtrsim N^{\max\{\frac{1}{2}, \frac{1}{p'}\}}$ for $1\le p<\infty$. 
To see that $\varphi_{u, \mcE}^p(N) \lesssim N^{\max\{\frac{1}{2}, \frac{1}{p'}\}}$, fix \(A \subset \mathbb{N}\) with \(\abs{A}= N \).\\
If  \( 2\le p \le \infty\) we use the Hausdorff--Young inequality  \eqref{hausdorff-young:1-2} and get that
 	\begin{equation*}
	\normvar{\sum_{n\in A} \varepsilon_n n^{-s}}_{\mathcal{H}_p}\le  \left(\sum_{n\in A} \abs{\varepsilon_n}^{p'} \right)^{\frac{1}{p'}} =N^{\frac{1}{p'}}. 
 	\end{equation*} 
For  \( 1\le p \le 2 \),  the inclusion \( \mathcal{H}_2 \subset \mathcal{H}_p \) gives 
 	\begin{equation*}
 		\normvar{\sum_{n\in A} \varepsilon_n n^{-s}}_{\mathcal{H}_p} \le  \normvar{\sum_{n\in A} \varepsilon_n n^{-s}}_{\mathcal{H}_2} =  N^{\frac{1}{2}}.
 	\end{equation*}

For the lower functions in (b)(i), 
combining \eqref{sum first N primes}  and \eqref{eq-norma-pot-de-2} we get  
 		\begin{equation*}
 			\varphi_{l,\mcE}^p(N) \lesssim 
 				N^{\min \left\{\frac 1 2, \frac 1 {p'} \right\}}.
 		\end{equation*} 
Similarly to the previous case, the reverse estimates follow from the Hausdorff--Young inequalities \eqref{hausdorff-young:2-infty}  and the inclusion $\mathcal{H}_p\subset \mathcal{H}_2$ (for $2\le p<\infty$). 

Finally, we deal with the case $p=\infty$. Note that from \eqref{eq-norma-dir-pol}  we have  $\varphi_{u,\mcE}^\infty (N) \ge \varphi_{u}^\infty (N) \ge \absito{\sum_{n=1}^N n^{i0}} = N$ and then the equality holds. For the lower functions,  \eqref{sum first N primes} shows  that $\varphi_{l}^\infty(N)=N$.

The inclusion $\mathcal{H}_\infty \subset \mathcal{H}_2 $ gives $\varphi_{l,\mcE}^\infty (N)\ge N^{\frac 1 2}$. For the reverse estimate, given \( N \in \mathbb{N}\) we  take \(k \) such that $2^{k-1}\le N \le 2^{k}$
and consider the $k$th Rudin--Shapiro polynomial (see Definition \ref{rud-shap})
\[ 
P_k(z)= \sum_{j=1}^{2^{k}} a_{kj}z^j, 
\]
with $a_{kj}=\pm1$. Since \( \normita{P_k}_{\infty} \le \sqrt{2^{k+1}} \), via the Bohr transform, we get 
\[
\varphi_{l,\mcE}^\infty(N) \le \sqrt{2^{k+1}} \le  2 \sqrt{N}, \]
which finishes the proof. 
\end{proof}

\section{Hardy spaces of general Dirichlet series} 
\label{sec-generales}

Fix \(\lambda := (\lambda_n) \) an unbounded sequence of strictly increasing non-negative real numbers, (called \emph{frequency}). A $\lambda$\emph{-Dirichlet} series is a formal sum of the form \( D = \sum_n a_n e^{-\lambda_n s},\) where \(a_n\) are complex coefficients and \(s\) is a complex variable. Any sum \( \sum_{n\in A}  a_n e^{-\lambda_n s} \) for a finite set $A\subset \NN$, is called $\lambda$-\emph{Dirichlet polynomial}.

As in \eqref{eq-norma-dir-pol}, for $1\le p<\infty$ we define the $p$-norm of a $\lambda$-Dirichlet polynomial $D$ with coefficients $(a_n)_{n\in A}$ as
\begin{equation}\label{eq-pnorm}
  \|D\|_{\mathcal{H}^\lambda_p}= \lim_{R \to \infty} \Big(\frac{1}{2R}\int_{-R}^R \Big|\sum_{n\in A} a_n e^{-\lambda_n it} \Big|^p dt\Big)^\frac{1}{p} \,.
\end{equation}
Now, the \emph{Hardy space $\mathcal{H}^{\lambda}_p$} is defined as the completion of the space of $\lambda$-Dirichlet polynomials  with the above norm. In the ordinary case, we  have an approach to Dirichlet series through  Fourier analysis in $\mathbb T^\infty$ and  Bohr's transform.  In the general case, there is also an approach through Fourier analysis on groups (see \cite[Definition~3.23 and Theorem~3.26]{DeSch}). We will not delve into details, but it's worth noting that this approach enables us to apply well-known results (such as Hölder, Littlewood or Hausdorff--Young inequalities, see \cite[Section~3.1]{CaDeMaSchSe2024}) to general Dirichlet series. 

The definition of $\mathcal{H}_\infty^\lambda$ makes use of the group approach to Hardy spaces. However, our main goal is to study properties of the basis $\mathcal B=\{e^{-\lambda_ns}\}$ that only requires computing the $\infty$-norm of Dirichlet polynomials. For $D = \sum_{n\in A} a_n e^{-\lambda_n s}$ with $A\subset \NN$ finite and scalars $(a_n)_{n\in A}$ we have 

\begin{equation}\label{norm-gral-infty}
\|D\|_{\mathcal{H}^\lambda_\infty} = \sup_{t \in \mathbb{R}} \Big| \sum_{n\in A} a_n e^{-i\lambda_n t} \Big|.    
\end{equation}

As in the ordinary case, we define $\mathcal A^\lambda_\infty$ as the completion of the space of Dirichlet polynomials with the above norm. 

Notice that, whenever $1\leq p < q \leq \infty$ we have the continuous incusions $\mcH_q^\lambda \subset \mcH_p^\lambda$, with $\Vert \cdot \Vert_{\mcH^\lambda_p} \leq \Vert \cdot \Vert_{\mcH^\lambda_q}$ that, contrary to what happens in the ordinary case,  may not be strict.

We begin by establishing the asymptotic behavior of the fundamental functions of the canonical basis $\mathcal B=\{e^{-\lambda_ns}\}$ of $\mathcal{H}^\lambda_p$, for a general $\lambda$. Recall that, for $1<p<\infty$, $\mathcal{B}$ is a Schauder basis of $\mathcal{H}_p^\lambda$ (see \cite[Theorem~4.16]{DeSch}).  Also, as happens in the ordinary case, $\mathcal B$ is a Markushevich basis if $p=1$ or $p=\infty$, in the latter case of $\mathcal A_\infty^\lambda$. To simplify the notation, we write  \begin{align*}
\varphi_{u}^{p,\lambda}(N):=\varphi_u[\mathcal B,\mathcal{H}_p^\lambda](N),\ \,&\qquad \varphi_{l}^{p,\lambda}(N):=\varphi_l[\mathcal B,\mathcal{H}_p^\lambda](N), \\ \varphi_{u,\mcE}^{p,\lambda}(N):=\varphi_{u,\mcE}[\mathcal B,\mathcal{H}_p^\lambda](N),\ \,&\qquad \varphi_{l,\mcE}^{p,\lambda}(N):=\varphi_{l,\mcE}[\mathcal B,\mathcal{H}_p^\lambda](N).
\end{align*}

The next result extends Theorem~\ref{teo-democ-ordinarias}. 

\begin{theorem}\label{teo-fi-series-generales}
Let \( \lambda \) be an arbitrary frequency. The fundamental functions of the canonical basis $\mcB=\{e^{-\lambda_n s}\}$ in $\mathcal{H}^{\lambda}_p$ satisfy:
 
\begin{enumerate}[\upshape (a)]
\item If $p=1$,\vspace{-10pt}
\begin{eqnarray*}
{\rm (i)}\quad  
\log(N) \lesssim \varphi_{l,\mcE}^{1,\lambda}(N) \le \varphi_{l}^{1,\lambda}(N) \lesssim N^{\frac{1}{2}},  
& \  & 
{\rm (ii)}\quad 
\varphi_{u}^{1,\lambda}(N)\approx \varphi_{u,\mcE}^{1,\lambda}(N) \approx N^{\frac{1}{2}}. 
\end{eqnarray*} \vspace{-20pt}
\item \label{eq-fiNcotas1}  If $1<p\le 2$, \vspace{-10pt}
\begin{eqnarray*} 
{\rm (i)}\quad 
N^{\frac{1}{p'}} \lesssim \varphi_{l,\mcE}^{p,\lambda}(N) \le \varphi_{l}^{p,\lambda}(N) \lesssim N^{\frac{1}{2}},
& \quad &
{\rm (ii)}\quad 
\varphi_{u}^{p,\lambda}(N) \approx \varphi_{u,\mcE}^{p,\lambda}(N) \approx N^{\frac{1}{2}}. 
\end{eqnarray*} \vspace{-20pt}
\item \label{eq-fiNcotas2} If $2 < p < \infty,$ \vspace{-10pt}
\begin{eqnarray*}   \hspace{7pt}%
{\rm (i)}\quad 
\varphi_{l,\mcE}^{p,\lambda}(N) \approx \varphi_{l}^{p,\lambda}(N) \approx N^{\frac{1}{2}}, \hspace{.7cm}
&  & 
{\rm (ii)}\quad 
N^{\frac{1}{2}} \lesssim \varphi_{u}^{p,\lambda}(N) \approx \varphi_{u,\mcE}^{p,\lambda}(N) \lesssim N^{\frac{1}{p'}}. 
\end{eqnarray*} \vspace{-20pt}
\item \label{eq-fiNcotas(d)}
If $p= \infty,$ \vspace{-10pt}
\begin{eqnarray*}
{\rm (i)}\quad  N^{\frac12} \lesssim \varphi_{l,\mcE}^{\infty,\lambda}(N) \le  \varphi_{l}^{\infty,\lambda}(N) =N &  &{\rm (ii)}\quad  
  \varphi_{u}^{\infty,\lambda}(N) \ = \ \varphi_{u,\mcE}^{\infty,\lambda}(N) \ =\  N. 
\end{eqnarray*}
\end{enumerate}
\end{theorem}

\begin{proof}
    We start by proving that for $1\le p < \infty$ the following holds 
    \begin{equation}\label{ineq1-fi-series-generales}
        \varphi_{l,\mcE}^{p,\lambda}(N) \le\varphi_{l}^{p,\lambda}(N) \lesssim N^{\frac{1}{2}} \lesssim \varphi_{u}^{p,\lambda}(N) \approx  \varphi_{u,\mcE}^{p,\lambda}(N).
    \end{equation}
	By definition of the lower and upper super-democracy functions, we only need to find sets  \( A \subset \mathbb{N} \) with \( \abs{A}=N \) satisfying 
	\begin{equation*}
		\normvar{\sum_{n\in A}  e^{-\lambda_n s}}_{\mathcal{H}_p^\lambda} \approx N^{\frac{1}{2}}.
	\end{equation*}
Since $(\lambda_n)$ is unbounded, we can take a subsequence $(\lambda_{n_k})_k$ such that $\lambda_{n_{k+1}}\ge 2 \lambda_{n_k}$ for all $k$. In this way, $(\lambda_{n_k})_k$ satisfies the hypothesis (of being a lacunary) required for the frecuency in \cite[Theorem~3.7.4]{grafakos2008classical} . Moreover, from the proof, we observe that the result in \cite{grafakos2008classical} remains true for arbitrary frequencies of the type. Then we have
$$ 
\normvar{\sum_{k=1}^{N}  e^{-\lambda_{n_k}s}}_{\mathcal{H}_p^\lambda} \approx N^{\frac 12}.
$$

On the other hand, for \( 1\le p \le 2\) the inclusion $\mathcal{H}_2^\lambda \subset \mathcal{H}_p^\lambda$ gives
\begin{equation}\label{ineq2-fi-series-generales}
    \varphi_{u,\mcE}^{p,\lambda}(N) \lesssim N^{\frac12}.
\end{equation}
While for $2\le p \le \infty$, the inclusion $\mathcal{H}_p^\lambda \subset \mathcal{H}_2^\lambda$ gives
\begin{equation}\label{ineq3-fi-series-generales}
    N^{\frac12} \lesssim \varphi_{l,\mcE}^{p,\lambda}(N).
\end{equation}
Combining \eqref{ineq1-fi-series-generales} and \eqref{ineq2-fi-series-generales} we obtain (a)(ii) and (b)(ii). And combining \eqref{ineq1-fi-series-generales} with \eqref{ineq3-fi-series-generales} we obtain~(c)(i).

From \cite[Theorem~1.1]{JamingKellaySaba2023} or \cite[Theorem~1]{hudson1992hardy} we get that for all  \(A \subset \NN\) with \( \abs{A}=N\)  
    \begin{equation*}
        \normvar{\sum_{n\in A} c_n e^{-\lambda_n s}}_{\mathcal{H}^\lambda_1} \ge C^{-1} \log (N).
    \end{equation*}
   whenever \( \abs{c_n}\ge 1 \) for all \(n\). This completes the proof of (a)(i).

The first estimate in (b)(i) follows using the Hausdorff--Young inequalities for general Dirichlet series (see, for example, \cite[Section~3.1]{CaDeMaSchSe2024}). The third estimate in (c)(ii) follows analogously. 
 
Finally, we deal with the case $p=\infty$. Note that, using the $\infty$-norm of Dirichlet polynomials, it is easy to see that $\varphi_{l}^{\infty,\lambda}(N)=\varphi_{u}^{\infty,\lambda} (N)=\varphi_{u,\mcE}^{\infty,\lambda} (N)=N$.  And the inclusion $\mathcal{H}_\infty \subset \mathcal{H}_2 $ gives $\varphi_{l,\mcE}^{\infty,\lambda} (N)\gtrsim N^{\frac 1 2}$.
 \end{proof}

We remark that inequalities of Theorem~\ref{teo-fi-series-generales} are sharp. Take, for example, $\lambda=(\log n)_n$ (i.e., ordinary Dirichlet series, see Theorem~\ref{teo-democ-ordinarias}) or $\lambda=(\log p_n)_n$ (the logs of the prime numbers), for which $\varphi_{u,\mcE}^{p,\lambda}(N)\approx N^{\frac 1 2 }$ for any finite $p$ and $\varphi_{l,\mcE}^{\infty,\lambda}(N)= N$, 
see \eqref{sum first N primes}. 

Now we see that, for $r>2$, any intermediate behavior of \ref{eq-fiNcotas2}(ii) is possible in $\mathcal{H}_r^\lambda$.

\begin{proposition}\label{prop-Bourgain}
For every $r>2$ and every $\frac{1}{2}<t<\frac{1}{r'}$ there exists a frequency $\lambda$ such that 
$$
\varphi_{u}^{r,\lambda}(N)\approx \varphi_{u,\mcE}^{r,\lambda}(N)\approx N^t.
$$
\end{proposition}
\begin{proof}
    Given $2<p<r$ (to be chosen later), from  the proof of \cite[Theorem~2]{Bo} we can take a frequency $\lambda \subset \mathbb{N}$ and constants $C_1$ and $C_2$ such that: 
\begin{itemize}
		\item for every choice of coefficients $(a_n) \subseteq \mathbb C$ and for every finite set $A\subset \NN$, we have
		\begin{equation} \label{freude}
		\Big\|\sum_{n\in A} a_n e^{-\lambda_n s} \Big\|_{\mathcal{H}_p^\lambda}
		\le C_1 \, \Big( \sum_{n \in A} \vert a_n \vert^{2} \Big)^{1/2};
		\end{equation}
		\item for each $j$ there exists $S_j\subset \lambda$ with $\abs{S_j}=\lfloor 4^{j/p} \rfloor$ such that
		\begin{equation} \label{freiheit}
		\Big\|\sum_{n\in S_j} e^{-\lambda_n s}\Big\|_{\mathcal{H}_r^\lambda}\ge C_2\, 2^{j(1/p-1/r)}\abs{S_j}^{1/2}.
		\end{equation}
 \end{itemize}
Actually, \cite[Theorem~2]{Bo} is stated for Fourier series, but since we are working with an integer frequency, for any $q$ and any finite sum, we have 
$$
\Big\|\sum_{n} a_n e^{-\lambda_n s} \Big\|_{\mathcal{H}_q^\lambda} = \Big\|\sum_{n} a_n z^{\lambda_n} \Big\|_{L_q(\mathbb T)}.
$$ 
Hence, the reformulation in terms of Dirichlet series is straightforward. 
From~\eqref{freiheit}, with $N=\lfloor 4^{j/p} \rfloor$, $j\in \mathbb N$, we get 
$$ 
\varphi_{u}^{r,\lambda}(N) \ge C_2 N^{{1}/{2}-{p}/{2r} } N^{1/ 2} = C_2 N^{1-{p}/{2r}}. 
$$ 
Take now an arbitrary $N\in \mathbb N$ and choose $j$ such that
$$
\lfloor 4^{j/p} \rfloor \le N \le \lfloor 4^{(j+1)/p} \rfloor.
$$
We then have
$$
\varphi_{u}^{r,\lambda}(N)\ge \varphi_{u}^{r,\lambda}({\lfloor 4^{j/p} \rfloor}) \ge C_2 \ \lfloor 4^{j/p} \rfloor^{1-{p}/{2r} } \ge \frac{C_2}{4^{1-p/2r}} \lfloor 4^{(j+1)/p} \rfloor^{1-{p}/{2r}}  \ge \frac{C_2}{4^{1-p/2r}} N^{1-{p}/{2r}},
$$
meaning that 
$$
\varphi_{u}^{r,\lambda}(N)\gtrsim N^{1-{p}/{2r}}.
$$ 
We choose $p=2r(1-t)$, so that the above exponent is $t$ (it is easy to check that $2<p<r$), so we have $\varphi_{u}^{r,\lambda}(N)\gtrsim N^t.$ 

To get the lower estimate, we take  $A\subset \mathbb N$ with $\abs{A}=N$. Using Littlewood's inequality \cite[Theorem~5.5.1]{Garling-book} (which applies to general Dirichlet series via the group approach, as mentioned above) and then \eqref{freude} we have
\begin{eqnarray*}
\Big\|{\sum_{n\in A} \varepsilon_n e^{-\lambda_{n} s} }\Big\|_{\mathcal{H}_r^\lambda} &\le&  \Big\|{\sum_{n\in A} \varepsilon_n e^{-\lambda_{n} s} }\Big\|_{\mathcal{H}_p^\lambda}^{p/r} \  \Big\|{\sum_{n\in A} \varepsilon_n e^{-\lambda_{n} s} }\Big\|_{\mathcal{H}_\infty^\lambda}^{1- {p/r}} \\ & \lesssim &  N^{p/(2r)} N^{1- {p/r}} =N^{1-{p}/{2r}}=N^t.  
\end{eqnarray*} This gives $\varphi_{u,\mcE}^{r,\lambda}(N)\lesssim N^t$ which, together with the previous inequality, finishes the proof.
\end{proof}

For ordinary Dirichlet series, the asymptotic behavior of the fundamental functions are completely determined by $p$. As shown in Theorem~\ref{teo-fi-series-generales} and Proposition~\ref{prop-Bourgain}, this is not the case for general Dirichlet series. In what follows, we present different families of frequencies  which allow us to give more precise estimates. 

\subsubsection*{$\mathbb Q$-linearly independent frequencies.} 

It is well-known that, if $\lambda$ is a $\mathbb Q$-linearly independent sequence, then $\mathcal{H}_p^\lambda = \mathcal{H}_2^\lambda$ for $1\le p<\infty$, with equivalent norms (see, for example, \cite[Theorem 3.10]{CaDeMaSch21}, where this is shown in the vector valued case). Also, $\mathcal{A}_\infty^\lambda = \ell_1$ with equivalent norms (identifying each Dirichlet series with the sequence of its coefficients), as shown in \cite[Theorem 4.7]{School20}. Just for completeness, we recall the argument for $p=\infty$ and Dirichlet polynomials.  

Take $\lambda$ a linearly independent frequency over $\mathbb Q$. By the Kronecker theorem, we know that for each $M\in \mathbb N$, the set $\{(e^{-i\lambda_{1}t},\dots, e^{-i\lambda_{M}t})\colon t\in \mathbb R\}$ is dense in $\mathbb T^M$. Therefore, taking $M$ such that $A\subset \{1,\dots,M\}$ we have 
\begin{equation*}
		\normvar{ \sum_{n \in A}  a_n e^{-\lambda_{n}s} }_{\mathcal{H}_\infty^\lambda} = \sup_{t\in \mathbb R} \Big|\sum_{n \in A}  a_n e^{-i\lambda_{n}t}\Big| = \sup_{z\in \mathbb T^M} \Big|\sum_{n \in A}  a_n z_n\Big| = \sum_{n\in A} \abs{a_n}.
\end{equation*}

\subsubsection*{Lacunary frequencies}
We say that a frequency $\lambda=(\lambda_n)$ is lacunary if there exists a constant $L>1$ such that $\lambda_{n+1} \ge L\lambda_{n}$ for all $n\in \mathbb{N}.$ Lacunary frequencies behave similarly to $\mathbb Q$-linearly independent frequencies. 
Indeed, repeating the steps of the proofs of Theorems~3.7.4 and~3.7.6 from~\cite{grafakos2008classical} (originally done for integer frequencies but also applicable to arbitrary real frequencies) yields the following.
\begin{theorem}\label{Theorem equivalence of h lambda for lacunary frecuencies}
	Let \( \lambda \) be a lacunary frequency.  Then, 
 \begin{align*}
    \mathcal{H}_p^\lambda & = \mathcal{H}_2^\lambda \ \ \text{for all } \ \ 1\le p < \infty \ \ \text{and}\\   
    \mathcal{A}_\infty^\lambda & =\ell_1,
 \end{align*}
with equivalent norms.
\end{theorem}

\subsubsection*{Frequencies containing arbitrarily long arithmetic progressions }
This case includes $\lambda_n=n$ (Fourier series), $\lambda_n=\log(n)$ (ordinary Dirichlet series) and also examples like the following:
\[ \lambda := \{ \log n : n = \mathfrak{p}^\alpha, \ \alpha_j \le j \}.\]
For these frequencies, considering the estimations obtained in Theorem~\ref{teo-fi-series-generales}, we see that the democracy and super-democracy functions behave just as for ordinary Dirichlet series. 

\begin{proposition} \label{prop-progressions}
Let $\lambda$ be a frequency containing arbitrarily  large  arithmetic progressions. Then, 
\begin{enumerate}[\upshape (a)]
    \item\label{artith:new-lower} $\varphi_{l}^{p,\lambda}(N)  \approx \varphi_{l,\mcE}^{p,\lambda}(N)  \approx	\begin{cases}
   \log(N) & \text{if}\ \ p = 1, \\
 	    N^{\frac{1}{p'}} & \text{if}\ \ 1<p\le 2,
 	\end{cases}$ \quad and \quad  $ \varphi_{l,\mcE}^{\infty,\lambda}(N)  \approx  N^{\frac{1}{2}}$.
	\item \label{artith:new-upper}
		$\varphi_{u}^{p,\lambda}(N)  \approx \varphi_{u,\mcE}^{p,\lambda}(N)  \approx	 N^{\frac{1}{p'}} \ \ \text{if}\ \  2 \le p < \infty$. 
\end{enumerate}
 
As a consequence, all democracy and super-democracy functions grow as in Theorem~\ref{teo-democ-ordinarias}.
\end{proposition}
\begin{proof}
     Assume first that $1\le p<\infty$. By Theorem~\ref{teo-fi-series-generales}, we only have to show that there are subsets \( A\subset \mathbb{N} \), with $|A|=N$, such that  
	\begin{equation}\label{eq-arbi-long-arit-progressions}
		\normvar{  \sum_{n \in A} e^{-\lambda_{n}s} }_{\mathcal{H}_p^\lambda} \approx \begin{cases}
			\log \abs{A} &\text{if}\ \ p=1,\\
			\abs{A}^{\frac{1}{p'}} &\text{if}\ \ p>1.
			\end{cases}
	\end{equation}
Given \( N \) we take \( A \subset \mathbb{N} \) with \( \abs{A} = N \) such that \( \{\lambda_n\colon n \in A\} \) is an arithmetic progression and write  \( A= \{ a+kb\colon k=1,\ldots,N \}\) for some \( a,b \in\mathbb{R} \). Then, 
	\begin{eqnarray*}
		\normvar{\sum_{n \in A} e^{-\lambda_{n}s} }_{\mathcal{H}_p^\lambda}^p  &= &\lim_{R\to\infty}\frac{1}{2R} \int_{-R}^{R}    \Big|\sum_{k=1}^{N} e^{-(a+kb)it} \Big|^p dt =
		\lim_{R\to\infty}\frac{1}{2R} \int_{-R}^{R}   \Big|\sum_{k=1}^{N} e^{-kb it}\Big|^p dt \\
    \\  &= & \lim_{R\to\infty}\frac{1}{2Rb} \int_{-Rb}^{Rb}    \Big|\sum_{k=1}^{N} e^{-k it} \Big|^p dt \approx   \normita{\mathbf{D}_N}_{L_p(\mathbb{T})}^p.
  \end{eqnarray*} 
Taking account \eqref{eq-norm-dir-ker} we see that \eqref{eq-arbi-long-arit-progressions} holds.\\
Thus, for $1\le p \le 2$, as stated in \ref{artith:new-lower}, we obtain the asymptotic behavior for the lower (super)-democracy functions. Also, for $2\le p<\infty$, we derive the estimates of the upper (super)-democracy functions, given in \ref{artith:new-upper}.
  
Now we deal with the case $p=\infty$, to complete \ref{artith:new-lower}. Given \(N \in \NN \), choose $k \in \NN$ such that $ 2^{k-1} \le N \le 2^k$ and consider the $k$th Rudin--Shapiro polynomial (see Definition \ref{rud-shap}):
\[ P_k(z)= \sum_{j=1}^{2^{k}} a_{kj}z^j, \]
where \(  a_{kj} =\pm 1\). We take   an arithmetic progression of length \( 2^k \), \( \{ a+jb\colon j=1,\ldots,2^k \} \subset \lambda\). Since  \( \normita{P_k}_{\infty} \le \sqrt{2^{k+1}}, \) we have
\begin{eqnarray*}
     \varphi_{l,\mcE}^{\infty,\lambda}(N) & \le& \normvar{\sum_{j=1}^{2^k} a_{kj} e^{-(a+jb)s}}_{\mcH_\infty^\lambda}  =
     \sup_{t\in \mathbb{R}} \Big|  \sum_{j=1}^{2^k} a_{kj} (e^{-bit})^j\Big|  \\
     &= & \sup_{z\in \mathbb{T}} \Big| \sum_{j=1}^{2^k} a_{kj} z^j\Big|     
 \le \sqrt{2^{k+1}} \le  2\sqrt{N}.          
\end{eqnarray*} 
Finally, from the behavior obtained in \ref{artith:new-lower} and \ref{artith:new-upper}, along with the rest of the estimates given in Theorem~\ref{teo-fi-series-generales} for arbitrary frequencies, it follows that the asymptotic behavior for this type of frequencies coincides with that given by ordinary bases.  
\end{proof}

\section{Unconditional and democratic properties of the canonical basis}
\label{sec-unconditional}

For ordinary Dirichlet series, the canonical basis $\mcB=\{n^{-s}\}$ is unconditional in $\mathcal H_p$ only for $p=2$. As mentioned in Remark \ref{rem-conseq}, $p=2$ is also the only case in which $\mathcal B$ is super-democratic and, furthermore,  $p=2$ and $p=\infty$ are the only cases in which $\mathcal B$ is democratic. 
For general Dirichlet series, there are situations where $\mathcal{H}_p^\lambda \approx \mathcal{H}_2^\lambda$ for $p\ne 2$. Thus, it is possible to have unconditionality and (super-)democracy for some $p\ne 2$. The natural question, then, is whether good properties of the basis in some $\mathcal{H}_p^\lambda$ necessarily imply that $\mathcal{H}_p^\lambda \approx \mathcal{H}_2^\lambda$. While, this is somewhat known for unconditionality, we present the standard argument here, as we have not found a direct reference. 

\begin{proposition}\label{prop-uncond}
Let $\lambda$ be a frequency.
The sequence \( \{e^{-\lambda_ns}\} \) is an  unconditional Schauder basis of \( \mathcal{H}_p^\lambda \) with $1\le p<\infty$ (respectively, of  $\mathcal A_{\infty}^\lambda$) if and only if \( \mathcal{H}_p^\lambda  = \mathcal{H}_2^\lambda \) (respectively,   $\mathcal A_{\infty}^\lambda = \ell_1$), with equivalent norms. 
\end{proposition}

\begin{proof} The \textit{if} part is obvious. For the \textit{only if} part, assume first that \( \{e^{-\lambda_ns}\} \) is unconditional of \(\mathcal{H}_p^\lambda \) for \(1\le p <\infty\), then  
\begin{equation}\label{eq-la misma de ayer}
    \normvar{\sum_{n=1}^{N}  a_n e^{-\lambda_{n} s }}_{\mathcal{H}_p^\lambda}  \approx \normvar{\sum_{n=1}^{N} w_n  a_n e^{-\lambda_{n} s }}_{\mathcal{H}_p^\lambda},
\end{equation} 
for all $N\in \mathbb N $, all $w=(w_1,\dots,w_N) \in \mathbb T^N$ and all $a_1,\dots,a_N \in \mathbb{C}$.
Averaging on $\mathbb T^N$ we have
			
\begin{eqnarray}
    \normvar{\sum_{n=1}^{N}  a_n e^{-\lambda_{n} s }}_{\mathcal{H}_p^\lambda} &  \approx & \int_{\mathbb T^N} \normvar{\sum_{n=1}^{N} w_n a_n e^{-\lambda_n s}      }^p_{\mathcal{H}_p^\lambda} dw 
    = \int_{\mathbb T^N} \lim_{R\to\infty} \frac{1}{2R} \int_{-R}^{R}  \Big|  \sum_{n=1}^{N} w_n a_n e^{-\lambda_n i t} \Big|^p dt \ dw  \nonumber\\
    &=& \lim_{R\to\infty} \frac{1}{2R} \int_{-R}^{R} 	\int_{\mathbb T^N} 
    \Big|\sum_{n=1}^{N} w_n a_n e^{-\lambda_n it} \Big|^p  dw \  dt \approx \lim_{R\to\infty}\frac{1}{2R} \int_{-R}^{R} dt \left(\sum_{n=1}^{N} \abs{a_n}^2 \right)^{\frac{p}{2}} \nonumber\\
    &=& \left(\sum_{n=1}^{N} \abs{a_n}^2 \right)^{\frac{p}{2}}, \label{eq-cuenta-kin}
\end{eqnarray}
where we first combined the dominated convergence and the Fubini--Tonelli theorems and then used the Khinchin--Steinhaus inequality \eqref{eq-kin}. 

Suppose now that \( \{e^{-\lambda_ns}\} \) is unconditional in \( \mathcal{A}_\infty^\lambda\). 
Then, equation \eqref{eq-la misma de ayer} holds for \(p=\infty.\) Taking, \(w_n=\overline{a}_n/\abs{a_n}\), we get 
    \begin{equation*}
        \normvar{\sum_{n=1}^{N}  a_n e^{-\lambda_{n} s }}_{\mathcal{H}_\infty^\lambda}  \approx \normvar{\sum_{n=1}^{N} \abs{a_n} e^{-\lambda_{n} s }}_{\mathcal{H}_\infty^\lambda} \ge \sum_{n=1}^N \abs{a_n}, 
    \end{equation*}
yielding that \( \{e^{-\lambda_ns}\} \) is equivalent to the canonical basis of $\ell_1$. 
\end{proof}

Regarding democracy, recall that for frequencies containing arbitrarily long arithmetic progressions, the fundamental functions behave as in the ordinary case. So, for this family of frequencies, we have democracy only for $p=2$. 
As already mentioned, we cannot expect this for arbitrary frequencies, but we may wonder if democracy in $\mathcal{H}_p^\lambda$ implies \( \mathcal{H}_p^\lambda  \approx \mathcal{H}_2^\lambda \).  What we can show is that this \emph{almost} occurs for even integer $p$, in the following sense. 

\begin{proposition}\label{prop-democ-aprox-h2}
Let $\lambda$ be a frequency.
    If  \( \mathcal B=\{e^{-\lambda_ns}\} \) is democratic for $\mathcal{H}^\lambda_{2k}$,  $k\in \mathbb N$, then 
    \begin{equation}\label{eq-democ-estim}
        \normvar{\sum_{n\in A}  a_n e^{-\lambda_{n} s }}_{\mathcal{H}_{2}^\lambda} \le \normvar{\sum_{n\in A}  a_n e^{-\lambda_{n} s }}_{\mathcal{H}_{2k}^\lambda}  \lesssim \log(|A|+1) \normvar{\sum_{n\in A}  a_n e^{-\lambda_{n} s }}_{\mathcal{H}_{2}^\lambda} 
    \end{equation} 
for all finite sets $A\subset \mathbb N$ and all scalars $(a_n)\subset \mathbb C$.
\end{proposition}

Before giving the proof of the above proposition, observe that from the estimates \eqref{eq-norm-dir-ker} for the norm of the Dirichlet kernel $\mathbf{D}_N$,  we have 
\begin{equation}\label{eq-comparing}
    \normvar{ \sum_{j=1}^{N}  (2^j)^{-s} }_{\mathcal{H}_{2k}}  \approx N^{\frac{1}{(2k)'}-\frac12} \normvar{ \sum_{j=1}^{N}  (2^j)^{-s} }_{\mathcal{H}_{2}} = N^{\frac{k-1}{2k}} \normvar{ \sum_{j=1}^{N}  (2^j)^{-s} }_{\mathcal{H}_{2}}.
\end{equation}   
Comparing \eqref{eq-democ-estim}  with \eqref{eq-comparing}  we see that, whenever \( \mathcal B=\{e^{-\lambda_ns}\} \) is democratic for $\mathcal{H}^\lambda_{2k}$,  the $2k$-norm is much closer to the $2$-norm than what happens in the ordinary case for $\{n^{-s}\}$ in $\mathcal{H}_{2k}$ (considering that in this case the basis is not democratic). 

\begin{proof}[Proof of Proposition~\ref{prop-democ-aprox-h2}]
   We recall the definition of the \emph{$k$th additive energy} of a set $L \subseteq \mathbb{R}$, given and denoted by
$$E_k(L)=|\{(\lambda,\mu)\in L^{2k}: \ \lambda_1+\cdots+\lambda_k=\mu_1+\cdots+\mu_k\}|.$$
In \cite[Section~3.3]{CaDeMaSchSe2024}, the following relationship between the $k$th additive energy and the $\mathcal{H}_{2k}^\lambda$-norm is established:
 \begin{equation}
     \Big\| \sum_{n: \lambda_n\in L} e^{-\lambda_n s} \Big\|_{\mathcal{H}^\lambda_{2k}}
=       E_k(L)^{\frac 1 {2k}}\,. \label{eq-Feli}
 \end{equation}   
Now, if  \( \{e^{-\lambda_ns}\} \) is democratic for $\mathcal{H}^\lambda_{2k}$,  by Theorem~\ref{teo-fi-series-generales} and \eqref{eq-Feli} we have  
$$ 
\abs{A}^{\frac{1}{2}} \approx \normvar{  \sum_{n \in A} e^{-\lambda_{n}s} }_{\mathcal{H}_{2k}^\lambda} =  E_k(\{\lambda_n\colon n\in A\})^{\frac 1 {2k}}\,,
$$ 
for any finite set $A\subset \mathbb N$.
This means that there are positive constants $C_1$ and $C_2$ such that $$ C_1\le \frac{ E_k(L)^{\frac 1 {2k}}}{\abs{L}^{\frac{1}{2}}} \le C_2$$ for all  finite sets $L\subset \lambda$. 
An application of \cite[Lemma~3.8]{CaDeMaSchSe2024} now gives the result. 
  \end{proof}

In the previous proof we used the fact that, by Theorem~\ref{teo-fi-series-generales}, \( \mcB=\{e^{-\lambda_n s}\} \) is democratic in $\mathcal{H}_p^\lambda$ if and only if both democracy functions behave like $N^{1/2}$. The following result says that the same happens if we change the assumption of democracy to SUCC. Moreover, it states that whenever $\mcB$ is SUCC,  it is (super-)democratic, which is not the usual case as \cite[Sec.~1.3, p 23.]{temlyakov2011greedy} shows.

\begin{theorem}\label{teo-qglc_dn} 
Let $\lambda$ be a frequency.
The basis \( \mcB=\{e^{-\lambda_n s}\} \) in $\mathcal{H}_p^\lambda$
is SUCC  if and only if it is super-democratic.\\ 
In addition, when this holds we have 
\begin{eqnarray} 
\varphi_{l,\mcE}^{p,\lambda}(N) \approx & N^{\frac{1}{2}} &\approx \varphi_{u,\mcE}^{p,\lambda}(N) \quad 1\le p<\infty \quad \text{or}  \label{succ=superdemo} \\
\varphi_{l,\mcE}^{\infty,\lambda}(N) \approx &N& = \varphi_{u,\mcE}^{\infty,\lambda}(N). \label{succ=superdemo_inf}
\end{eqnarray} 
\end{theorem}

\begin{proof} For $1\le p <\infty$, we only have to prove that \( \{e^{-\lambda_n s}\} \) is SUCC if and only if it satisfies the estimations in \eqref{succ=superdemo}. The if part is clear. Conversely, suppose that \( \{e^{-\lambda_n s}\} \) is $\bC$-SUCC. Take  $A \subset \mathbb{N}$  with $\abs{A}=N$ and $\varepsilon \in \mcE_A$ such that
\begin{equation}\label{norm greater than KD(N)}
  	 \normvar{\sum_{n\in A} \varepsilon_n e^{-\lambda_{n} s}}_{\mathcal{H}_p^\lambda} \geq \frac{1}{2} \varphi_{u,\mcE}^{p,\lambda}(N).
  \end{equation}
By the Khinchin inequality \eqref{eq-kinchin} and proceeding as in \eqref{eq-cuenta-kin}, we have
\begin{equation*}
  \mathbb{E}\normvar{\sum_{n \in A} \pm \varepsilon_{n} e^{-\lambda_{n} s}}_{\mathcal{H}_{p}^\lambda} \le \mathbf{K}_p|N|^{\frac{1}{2}},
\end{equation*}
so there must exist signs \( \delta_n = \pm 1 \) \((n \in A)\) such that 
\begin{equation}\label{sum delta signs A lower sqrN}
  \normvar{\sum_{n \in A} \delta_n \varepsilon_{n} e^{-\lambda_{n} s}}_{\mathcal{H}_{p}^\lambda} \le \mathbf{K}_p |N|^{\frac{1}{2}}.
\end{equation}
Now define \( A_1:= \{n\colon \delta_n=1\} \) and \( A_2:= \{n\colon \delta_n=-1\}\), which form a disjoint partition of $A$. By \eqref{norm greater than KD(N)} and the triangle inequality, we have 
 \begin{equation}\label{norm Aj ge K/2 D(N)}
  	 \normvar{ \sum_{n \in A_1} \varepsilon_{n} e^{-\lambda_{n} s}}_{\mathcal{H}_p^\lambda} \ge \frac{1}{4} \varphi_{u,\mcE}^{p,\lambda}(N) \quad \text{ or } \quad \normvar{\sum_{n \in A_2} \varepsilon_{n} e^{-\lambda_{n} s}}_{\mathcal{H}_p^\lambda} \ge \frac{1}{4} \varphi_{u,\mcE}^{p,\lambda}(N).
\end{equation}
Suppose that the first inequality holds (the other case is analogous). As  \( \{e^{-\lambda_n s}\} \) is $\bC$-SUCC, 
\begin{align*}
   \normvar{ \sum_{n \in A_1} \varepsilon_{n} e^{-\lambda_{n} s}}_{\mathcal{H}_p^\lambda}  \le  \bC \normvar{\sum_{n \in A} \delta_n \varepsilon_{n} e^{-\lambda_{n} s}}_{\mathcal{H}_{p}^\lambda} \nonumber
\end{align*}
The above inequality with \eqref{sum delta signs A lower sqrN} and \eqref{norm Aj ge K/2 D(N)} give 
\begin{equation*}
\varphi_{u,\mcE}^{p,\lambda}(N) \le 4\bC \mathbf{K}_p N^{\frac{1}{2}}.
\end{equation*}
Combining this with Theorem~\ref{teo-fi-series-generales}, we obtain \( \varphi_{u,\mcE}^{p,\lambda}(N)\approx N^{\frac{1}{2}}. \)
  
Let us see the estimate for the lower function, \( \varphi_{l,\mcE}^{p,\lambda}(N) \approx N^{\frac{1}{2}}. \) Take $A \subset \mathbb{N}$ with $\abs{A}\ge N$ and \(\varepsilon \in \mcE_A\) such that 
\begin{equation}\label{norm lower than Kd(N)}
  \normvar{\sum_{n\in A} \varepsilon_n e^{-\lambda_{n} s}}_{\mathcal{H}_p^\lambda} \leq  2 \varphi_{l,\mcE}^{p,\lambda}(N).
\end{equation}
Again, thanks to the Khinchin inequality \eqref{eq-kinchin}  we can find signs \( \delta_n = \pm 1 \) such that,  
\begin{equation}\label{sum delta signs A greater sqrN}
  	\normvar{\sum_{n \in A} \delta_n \varepsilon_{n} e^{-\lambda_{n} s}}_{\mathcal{H}_{p}^\lambda} \ge \mathbf{k}_p^{-1} |N|^{\frac{1}{2}}.
  \end{equation}
    As before, taking the disjoint partition of $A$, \( A_1:= \{n\colon \delta_n=1\} \) and \( A_2:= \{ n\colon \delta_n=-1\}\), we see that one of the following occurs:
\begin{equation}\label{norm Aj ge sqrN 2}
  \normvar{ \sum_{n \in A_1} \varepsilon_{n} e^{-\lambda_{n} s}}_{\mathcal{H}_p^\lambda} \ge \frac{\mathbf{k}_p^{-1}}{2} N^\frac{1}{2} \quad \text{ or } \quad \normvar{\sum_{n \in A_2} \varepsilon_{n} e^{-\lambda_{n} s}}_{\mathcal{H}_p^\lambda} \ge \frac{\mathbf{k}_p^{-1}}{2} N^{\frac{1}{2}}.
\end{equation}
  Now, by the $\bC$-SUCC property, 
\begin{align*}
  \normvar{ \sum_{n \in A_1} \varepsilon_{n} e^{-\lambda_{n} s}}_{\mathcal{H}_p^\lambda}  \le  \bC \normvar{\sum_{n \in A} \varepsilon_{n} e^{-\lambda_{n} s}}_{\mathcal{H}_{p}^\lambda} 
  \end{align*}
and, assuming that the first inequality of \eqref{norm Aj ge sqrN 2} holds, a combination of \eqref{norm lower than Kd(N)} and \eqref{sum delta signs A greater sqrN} yields 
  \begin{equation*}
  	\frac{\mathbf{k}_p^{-1}}{4\bC} N^{\frac{1}{2}} \le  \varphi_{l,\mcE}^{p,\lambda}(N).
  \end{equation*}
  This together with Theorem~\ref{teo-fi-series-generales} give the desired estimate.

For \( p = \infty\), Theorem~\ref{teo-fi-series-generales}\ref{eq-fiNcotas(d)} shows that $\{e^{-\lambda_{n} s}\}$ is always democratic with \(\varphi_{u}^{\infty,\lambda}(N)=\varphi_{l}^{\infty,\lambda}(N)=N\). As a basis is democratic and SUCC if and only if it is a super-democratic \cite[Proposition~5.1]{albiac2021dissertationes}, the conclusion follows. As a consequence, we obtain that $\varphi_{l,\mcE}^{\infty,\lambda}(N) \approx \varphi_{u,\mcE}^{\infty,\lambda}(N) =N$, and the proof is complete.  
\end{proof}

\section{Final remarks}
\label{sec-greedy}

In this section, we show that the basis of the ordinary Dirichlet series $\mcB=\{n^{-s}\}$ lacks the most commonly studied greedy properties whenever $p\ne 2$. In fact, we prove this for a wider class of Dirichlet series. To avoid lengthening this section by including the details of the definitions, we give a brief idea of the notions we study and establish how they relate. Then we apply the results obtained in the previous sections. The omitted details can be found in ~\cite{temlyakov1998greedy, temlyakov2011greedy}, the recent monograph~\cite{albiac2021dissertationes}, and the references therein.  

Given a basis $\mcB=(\xx_n)_{n}$ and $x\in \XX$, a finite set $A\subset \NN$ is a greedy set of $x$ if for all $j\in A$ and all $k\in \NN\setminus A$
\begin{align*}
|\xx_j^*(x)|\ge& |\xx_k^*(x)|, 
\end{align*}
where $\mcB^*=(\xx^*_n)_{n}$ is the dual basis of $\mcB$. Loosely speaking, the basis is \emph{greedy} if for every greedy set $A$ of $x$, $P_A(x)=\sum_{j\in A} \xx^*_j(x)\xx_j$ is the best approximation to $x$ through vectors whose supports have at most the cardinal of $A$. The basis is
\emph{quasi-gredy} if $\|P_A\|$ is uniformly bounded. This is equivalent to saying that the projections converge to the vector when the cardinal of the greedy sets tends to infinity, even though $P_A$ is neither linear nor continuous for any $|A|\ge 1$ (see e.g.,  \cite[Theorem~1.29]{temlyakov2011greedy}). 
A basis is greedy if and only if it is Schauder, unconditional, and democratic~\cite{konyagin1999remark}. This is the case, for example, of the canonical basis of $\ell_p$ for $1\le p<\infty$, or the normalized Haar basis of $L_p[0, 1]$ for $1<p<\infty$. Notice that every unconditional basis is quasi-greedy while there are quasi-greedy bases (even enjoying the super-democracy property) that are conditional, as shown in~\cite{konyagin1999remark}, see also \cite[p.~35]{temlyakov2011greedy}.

Intermediate structures arise from the works of Dilworth, Kalton, Kutzarova and Temlyakov under the name of (Schauder) \emph{almost greedy} bases, characterized as those that are quasi-greedy and democratic 
\cite{dilworth2003thresholding} and \emph{semi-greedy} bases~\cite{dilworth2003existence} given in terms of Chevishev-greedy type algorithms, which arise naturally in the setting. By the nature of the algorithm that uses approximations by sums of $m$-terms, Markushevich bases provide to the context a more flexible structure than Schauder bases do. 

For Schauder basis, the proof of the equivalence between semi-greedy and almost greedy basis was carried out in two steps, \cite[Theorem~3.2 and Theorem 3.6]{dilworth2003existence} and \cite[Theorem~1.10]{berna2019equivalence}. To include the cases $p=1$ and $p=\infty$ we need the equivalence for Markushevich basis \cite[Corollary~4.3]{berasategui2023weak}. 

Concerning quasi-greedy bases, a combination of \cite[Theorem~4.13]{albiac2021dissertationes} and \cite[Theorem~3.4]{albiac2022studia} yields that they are \emph{nearly unconditional} bases, a notion introduced by Elton in~\cite{elton1978weakly}. Recently, in \cite[Theorem~2.6]{albiac2023elton}, it was shown that the latter are equivalent to \emph{quasi-greedy bases for largest coefficients}, that in particular are SUCC.    

The diagram below summarizes the relationships described above.
\vskip -.2cm 
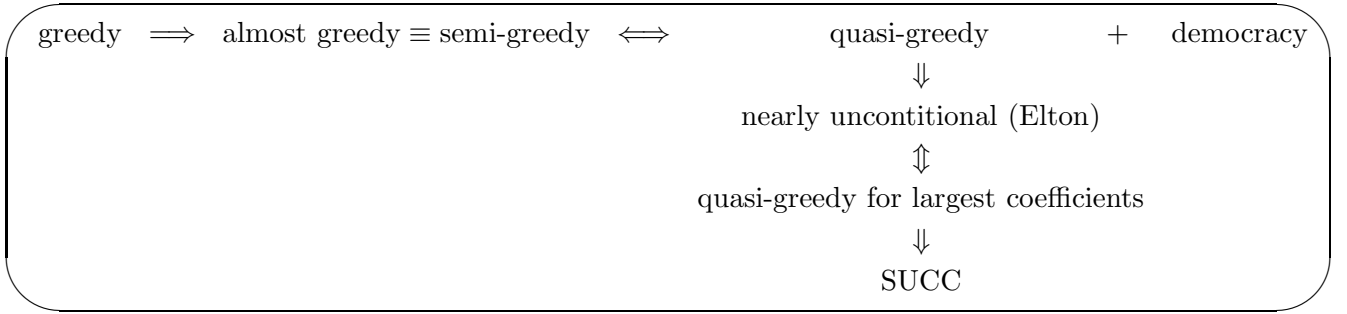
\begin{center}
\ovalbox{
$
\begin{array}{cccccc}
\text{greedy} & \Longrightarrow & \text{almost greedy} \equiv \text{semi-greedy} & \Longleftrightarrow & \qquad \qquad \text{quasi-greedy} \qquad\qquad  +&  \text{democracy}\\
& & & & \Downarrow &  \\ 
& & & & \text{nearly uncontitional (Elton)} &  \\
& & & & \Updownarrow &  \\ 
& & & & \text{quasi-greedy for largest coefficients} &\\
& & & & \Downarrow &  \\ 
& & & & \text{SUCC}
\end{array}
$}
\end{center}
\medskip

As already noticed, a basis is (super-)democratic if and only if the associated upper and lower (super-)democracy functions are asymptotically equivalent. Thus, our previous results yield the following consequences for the greedy properties mentioned above for $\mcB=\{e^{-\lambda_n s}\}$ in $\mcH^\lambda_p$ when $\lambda=(\lambda_n)$ is a frequency containing arbitrarily large arithmetic progression. 

\begin{proposition} \label{Prop:succ sii super-demo}
Let $\lambda=(\lambda_n)$ be a frequency containing arbitrarily large arithmetic progression. Then,  
$\mcB=\{e^{-\lambda_n s}\}$ in $\mcH^\lambda_p$ is SUCC if and only if $p=2$.
\end{proposition}

\begin{proof}
Suppose that $\mcB$ is SUCC. Then, by Theorem~\ref{teo-qglc_dn}, it is super-democratic and then either \eqref{succ=superdemo} or \eqref{succ=superdemo_inf} holds. Therefore, Proposition \ref{prop-progressions} and Theorem~\ref{teo-democ-ordinarias} imply that $p=2$. 
\end{proof}

We end this section summarizing the situation in the following Corollary.

\begin{corollary} \label{cor-supercoro}
Let $\lambda=(\lambda_n)$ be a frequency containing arbitrarily large arithmetic progression. Then, $\mcB=\{e^{-\lambda_n s}\}$ in $\mcH^\lambda_p$ satisfies any of the properties {\rm (i) --  (vii)} below if and only if $p=2$, with
\begin{enumerate}[\upshape (i)]
    \item $\mcB$ is greedy.
    \item $\mcB$ is almost greedy.
    \item $\mcB$ is semi-greedy.
    \item $\mcB$ is quasi-greedy.
    \item $\mcB$ is nearly unconditional.
    \item $\mcB$ is quasi-greedy for largest coefficients.
    \item $\mcB$ is SUCC.
\end{enumerate}
\end{corollary}

Notice that, in particular, Proposition~\ref{Prop:succ sii super-demo} and Corollary~\ref{cor-supercoro} hold for $\mcB=\{n^{- s}\}$ in $\mcH_p$ and $\mcB=\{e^{int}\}_{n\in \mathbb Z}$ in $L_p(\mathbb T)$. The last case was extensively studied, from the origins of the theory to the present day, see among many others \cite{temlyakov1998greedy} and \cite{berna2017lebesgue}.

\bibliographystyle{abbrv}
\bibliography{main}

\end{document}